\documentclass{amsart}                   
\usepackage{amsfonts}
\usepackage{amsmath}
\usepackage{amscd}
\usepackage{amssymb}
\usepackage[dvips]{graphicx}
\usepackage[dvips]{epsfig}
\usepackage{color}
\usepackage{hyperref}
\usepackage{url}

\def\Fraisse{Fra\"{\i}ss\'e\ }

\def\mb#1{\mathbb{#1}}
\def\mc#1{\mathcal{#1}}

\def\pigeonhole#1#2#3#4{\mb{#3} \ARROW{} (\mb{#2})^{\mb{#1}}_{#4}}

\def\sub#1{_{_{#1}}}
\def\sup#1{^{^{#1}}}

\def\age#1{\text{\rm Age}(\mb{#1})}
\def\flim{\text{\rm FLim}}

\def\Lim#1{\underset{#1}{\text{\rm Lim}}\ }

\def\R#1{\mb{R}\sup{#1}}                      
\def\N{\mb{N}}

\def\Comb#1#2#3{\left(\begin{array}{c} {#1} \\ {#2}\end{array}\right)\sub{#3}}
\def\comb#1#2{\Comb{\mb{#1}}{\mb{#2}}{}}

\def\aut#1{\text{\rm Aut}(\mb{#1})}

\def\ARROW#1{
\text{\begin{picture}(35,18)(15,15)         
            \put(30,22){$_{#1}$}
            \put(18,17){\vector(1,0){30}}
\end{picture}}}

\newcounter{numero}

\newcounter{letra}
\newcommand{\Letra}{\medskip \setcounter{letra}{1}(\alph{letra}) }
\newcommand{\letra}{\medskip \addtocounter{letra}{1}(\alph{letra}) }

\newcounter{romnumero}
\newcommand{\Romnumero}{\setcounter{romnumero}{1}(\roman{romnumero}) }
\newcommand{\romnumero}{\addtocounter{romnumero}{1}(\roman{romnumero}) }

\newcounter{bibnumero}

\newtheorem{teo}{Theorem}[subsection]                  
\newtheorem{lema}[teo]{Lemma}
\newtheorem{prop}[teo]{Proposition}
\newtheorem{cor}[teo]{Corollary}
\newtheorem{quest}[teo]{Question}
\newtheorem{quests}[teo]{Questions}

\def\bteo{\begin{teo}}
\def\eteo{\end{teo}}
\def\bprop{\begin{prop}}
\def\eprop{\end{prop}}
\def\bcor{\begin{cor}}
\def\ecor{\end{cor}}
\def\blema{\begin{lema}}
\def\elema{\end{lema}}
\def\bquest{\begin{quest}}
\def\equest{\end{quest}}
\def\bquests{\begin{quests}}
\def\equests{\end{quests}}

\theoremstyle{definition}                           
\newtheorem{definition}[teo]{Definition}

\newtheorem{definitions}[teo]{Definitions}
\newtheorem{ejems}[teo]{Examples}
\newtheorem{ejem}[teo]{Example}

\newtheorem{rem}[teo]{Remark}

\def\bdeff{\begin{definition}\rm }
\def\edeff{\end{definition}}
\def\bdefs{\begin{definitions}\rm }
\def\edefs{\end{definitions}}
\def\bejem{\begin{ejem}\rm }
\def\eejem{\end{ejem}}
\def\bejems{\begin{ejems}\rm }
\def\eejems{\end{ejems}}

\newenvironment{dem}{{\it Proof:} \rm\hskip3mm }{\hfill$\square$\vskip5mm}       
\newenvironment{sketch}{{\it Sketch of the Proof:}\rm\hskip3mm }{\hfill$\square$\vskip5mm}

\def\bdem{\begin{dem}}
\def\edem{\end{dem}}

\def\bsketch{\begin{sketch}}
\def\esketch{\end{sketch}}

\usepackage[english]{babel}
\usepackage{blindtext,tikz}
\usetikzlibrary{calc}

\begin{document}



\author{Jos\'e G.  Mijares}\thanks{For this research; Jos\'e G. Mijares was partially supported 
by the  project \textit{Din\'amica Topol\'ogica en Espacios Uniformes} ID PPTA 00005027, ID PRY 004863, 
Pontificia Universidad Javeriana. Bogot\'a, Colombia. Gabriel Padilla was partially supported by the project
\textit{Estructuras topol\'ogicas} ID 17426, Universidad Nacional de Colombia. Bogot\'a, Colombia.}
\address{Department of Mathematics, University of Denver. 2280  S Vine St., Denver, CO 80208, USA}
\email{Jose.MijaresPalacios@du.edu}

\author{Gabriel Padilla}
\address{Departamento de Matem\'aticas, Universidad Nacional de Colombia, K30-Cl45, 
Edif. 404, Ofic. 315-404. Bogot\'a 101001000 Colombia. (+571)3165000 ext 13166}
\email{gipadillal@unal.edu.co}

\title{A Ramsey space of infinite polyhedra\\ and the random polyhedron}



\begin{abstract}

In this paper we introduce a new topological Ramsey space $\mathcal P$ whose elements are infinite ordered 
polyhedra. The corresponding familiy $\mathcal{AP}$ of finite approximations can be viewed as a class of 
finite structures. It turns out that the closure of $\mathcal{AP}$ under isomorphisms is the class 
$\mathcal{KP}$ of finite ordered polyhedra. Following \cite{nerod3}, we show that $\mathcal{KP}$ is a Ramsey 
class. Then, we prove a universal property for ultrahomegeneous polyhedra and introduce the 
\textit{(ordered) random polyhedron}, and prove that it is the \Fraisse limit of $\mathcal{KP}$; hence the 
group of automorphisms of the ordered random polyhedron is extremely amenable (this fact is deduced from 
results of \cite{kpt}). Later,  we present a countably infinite family of topological Ramsey subspaces of 
$\mathcal P$; each one determines a class of finite ordered structures which turns out to be a Ramsey class. 
One of these subspaces is Ellentuck's space; another one is associated to the class of finite ordered graphs 
whose \Fraisse limit is the random graph.  The  \Fraisse limits of these classes are not 
pairwise isomorphic as countable structures and none of them is isomorphic to the random polyhedron. Finally, 
following \cite{kpt}, we calculate the universal minimal flow of the (non ordered) random polyhedron as well 
as the universal minimal flows of the (non ordered) random structures associated to our family of topological 
Ramsey subspaces of $\mathcal P$.
\end{abstract}
\maketitle


\section{Introduction}

A polyhedron is a geometric object built up through a finite or countable number of suitable amalgamations of convex hulls of finite sets; polyhedra are generated in this way by simplexes. Simplicial morphisms are locally linear maps that preserve vertices. An ordered polyhedron is a polyhedron for which we have imposed a linear order on the set of its vertices. As we only consider order-preserving morphisms, ordered polyedra are rigid, i. e., admit no non-trivial automorphisms; this is an easy consequence of the well order principle. In this paper we define a new topological Ramsey space 
(see \cite{todo}) whose elements are essentially infinite ordered polyhedra.

\medskip

The theory of topological Ramsey spaces is developed in \cite{todo}, and was pioneered by the work \cite{ellen} of Ellentuck's. In Section \ref{def ramsey spaces} we will describe the fundamental concepts of that theory. In Section \ref{def ramsey space polyhedra} we will define our new topological Ramsey space $\mathcal P$. The closure under isomorphisms of the corresponding family $\mathcal{AP}$ of finite approximations (viewed as a class of 
finite structures) turns out to be the class 
$\mathcal{KP}$ of finite ordered polyhedra.
In Section \ref{finite pol Ramsey}, following \cite{nerod3}, we prove the Ramsey property for the class 
$\mathcal{KP}$.  We also prove a universal property for ultrahomegeneous polyhedra and show that the automorphism 
group of  the \Fraisse limit of $\mathcal{KP}$ is extremely amenable, following \cite{kpt}. A description of this \Fraisse limit is given in Section  \ref{ssection infinite random polyhedron}; we call it the \textit{ordered random polyhedron}. In Section \ref{subespacios}, we introduce a countable family  
$\{\mathcal P(k)\}_{ k>0}$ of topological Ramsey subspaces of  $\mathcal P$. Each $\mathcal P(k)$ determines 
a class $\mathcal {KP}(k)$ of finite ordered structures which turns out to be a Ramsey class. The automorphism
group of its  \Fraisse limit 
is therefore extremely amenable. For instance, $\mathcal P(1)$  coincides with Ellentuck's space 
(see the definition below). The corresponding Ramsey class is off course the class of finite linearly ordered sets
whose  \Fraisse limit is $(\mathbb Q, \leq)$. On the other hand,  the Ramsey class associated to
$\mathcal P(2)$ is the class of finite ordered graphs whose  \Fraisse limit is the ordered random graph. 
It is worth mentioning that the  \Fraisse limits of the classes $\mathcal {KP}(k),\ k>0$, are not 
pairwise  isomorphic as countable structures, and none of them is isomorphic to the ordered random polyhedron.
Finally, following \cite{kpt}, we calculate the universal minimal flow of the (non ordered) random polyhedron 
as well as the universal minimal flows of the (non ordered) random structures associated to our family of 
topological Ramsey subspaces of $\mathcal P$.

\medskip

In brief, we introduce some new topological Ramsey spaces 
associated to polyhedra and related geometric and combinatorial objects, and study their relation with Ramsey classes of finite 
(ordered) structures and the auto\-morphism groups of their \Fraisse limits, and their universal minimal flows. 

\subsection*{Notation}

Given a countable set $A$, we will adopt the following notation throughout the paper. 
Let $A$ be a countable set and  $X\subseteq A$; then $|X|$ denotes the cardinality of $X$ and:
\begin{itemize}
\item $A^{[k]} = \{X\subseteq A : |X|=k\}$, for every $k\in\mathbb N$.
\item $A^{[\leq k]} = \{X\subseteq A : |X|\leq k\}$, for every $k\in\mathbb N$.
\item $A^{[<\infty]} = \{X\subseteq A : |X|<\infty\}$.
\item $A^{[\infty]} = \{X\subseteq A : |X|= \infty\}$.
\end{itemize}

\section{Ramsey spaces}\label{def ramsey spaces}
The definitions and results throughout this section are
taken from \cite{todo}.  

\subsection{Metrically closed spaces and approximations}\label{ MC spaces} 
Consider a triplet of the form $(\mathcal{R}, \leq,
r )$, where $\mathcal{R}$ is a set, $\leq$
is a quasi order on $\mathcal{R}$ and $r: \mathbb{N}\times\mathcal{R}\rightarrow \mathcal{AR}$
 is a function with range $\mathcal{AR}$. For every $n\in \mathbb{N}$ and every $A\in \mathcal{R}$, 
 let us write 
\begin{equation}\label{eq axiomatic approximations}
      r_n(A) := r(n,A)
\end{equation}
We say that $r_{n}(A)$ is {\bf the} $n$th {\bf approximation of} $A$. 
We will reserve capital letters $A,B\dots$ for elements in $\mc{R}$ while lowercase letters $a,b\dots$
will denote elements of $\mc{AR}$. 

\medskip
 
In order to capture the combinatorial structure
 required to ensure the provability of an Ellentuck type Theorem, some assumptions on
 $(\mathcal{R}, \leq, r)$ will be imposed. The first is the following:\\[2mm]
\noindent \textbf{(A.1)}
\begin{itemize}
	\item[{(A.1.1)}]For any $A\in \mathcal{R}$, $r_{0}(A) = \emptyset$.
	\item[{(A.1.2)}]For any $A,B\in \mathcal{R}$, if $A\neq B$ then
	$(\exists n)\ (r_{n}(A)\neq r_{n}(B))$.
	\item[{(A.1.3)}]If $r_{n}(A) =r_{m}(B)$ then $n = m$ and $(\forall i<n)\ (r_{i}(A) = r_{i}(B))$.
\end{itemize}
Take the discrete topology on $\mathcal{AR}$ and endow 
$\mathcal{AR}^{\mathbb{N}}$ with the product topology; this is the metric space of all the sequences of 
elements of $\mathcal{AR}$. The set $\mathcal{R}$ can be identified with the corresponding image 
in $\mathcal{AR}^{\mathbb{N}}$. We will say that $\mathcal{R}$ is
{\bf metrically closed} if, as a subspace $\mathcal{AR}^{\mathbb{N}}$ with the
inherited topology, it is closed. The basic open sets generating the metric topology on $\mathcal{R}$ inherited from the product topology of $\mathcal{AR}^{\mathbb{N}}$ are of the form:
\begin{equation}\label{eq basic opens in ARN}
		[a] = \{B\in \mathcal{R} : (\exists n)(a =
	r_{n}(B))\}
\end{equation}
where $a\in\mathcal{AR}$.

\medskip

 Let us define the {\bf length} of $a$,
as the unique integer $|a|=n$ such that $a = r_{n}(A)$ for some $A\in
\mathcal{R}$. For every $n\in\mathbb N$, let 
\begin{equation}
    \mathcal{AR}_n : = \{a\in\mathcal{AR} : |a|=n\} 
\end{equation}
Hence, 
\begin{equation}
    \mathcal{AR} = \bigcup_{n\in\mathbb N} \mathcal{AR}_n
\end{equation}

The {\bf Ellentuck type neighborhoods} are of the form:
\begin{equation}\label{eq basic Ellentuck opens}
		[a,A] = \{B\in [a] : B\leq A\}=
		\{B\in \mc{R} : (\exists n)\ a=r_n(B) \ \&\ B\leq A\}
\end{equation}
where $a\in\mathcal{AR}$
and $A\in \mathcal{R}$. 

\medskip

We will use the symbol $[n,A]$ to abbreviate $[r_{n}(A),A]$. 

\medskip

Let 
\begin{equation}
    \mathcal{AR}(A) = \{a\in\mathcal{AR} : [a,A]\neq\emptyset\}
\end{equation}
Given a neighborhood $[a,A]$ and $n\geq |a|$, let $r_n[a,A]$ be the image of $[a,A]$ by the function $r_n$, i.e., 
\begin{equation}\label{eq image of basic Ellentuck opens}
		r_n[a,A] =
		\{r_n(B) : B\in [a,A] \}
\end{equation}

\subsection{Ramsey sets}\label{ssection Ramsey sets}
A set $\mathcal{X}\subseteq \mathcal{R}$
is \textbf{Ramsey} if for every neighborhood $[a,A]\neq\emptyset$
there exists  $B\in [a,A]$ such that $[a,B]\subseteq \mathcal{X}$
or $[a,B]\cap \mathcal{X} = \emptyset$. A set
$\mathcal{X}\subseteq \mathcal{R}$ is \textbf{Ramsey null} if for
every neighborhood $[a,A]$ there exists  $B\in [a,A]$ such that
$[a,B]\cap \mathcal{X} = \emptyset$.

\subsection{Topological Ramsey spaces}\label{ssection Ramsey spaces}
We say that $(\mathcal{R}, \leq,r)$ is a \textbf{topological\linebreak Ramsey space} iff subsets of
$\mathcal{R}$ with the Baire property are Ramsey and\linebreak meager subsets of $\mathcal{R}$ are Ramsey null.

\medskip

Given $a,b\in\mathcal{AR}$, write 

\begin{equation}
a\sqsubseteq b\ \mbox{ iff } (\exists A\in\mathcal R)\ (\exists m,n\in\mathbb N)\  m\leq n, a=r_m(A)\mbox{ and } b= r_n(A).
\end{equation}

\medskip

By A.1, $\sqsubseteq$ can be proven to be a partial order on $\mathcal{AR}$.

\medskip

\noindent \textbf{(A.2)}\ [{\bf Finitization}] There is a quasi order $\leq_{fin}$ on
$\mathcal{AR}$ such that:
\begin{itemize}
    \item[(A.2.1)] $A\leq B$ iff
    $(\forall n)\ (\exists m) \ \ (r_{n}(A)\leq_{fin} r_{m}(B))$.
    \item[(A.2.2)] $\{b\in \mathcal{AR} : b\leq_{fin} a\}$ is finite, for every
    $a\in \mathcal{AR}$.
\item[(A.2.3)] If $a\leq_{fin} b$ and $c \sqsubseteq a$ then there is $d \sqsubseteq b$ such that $c \leq_{fin} d$.
\end{itemize}

\medskip

Given $A\in\mathcal{R}$ and $a\in\mathcal{AR}(A)$,
we define the {\bf depth} of $a$ in $A$ as 
\begin{equation}\label{eq depth of a segment}
		\operatorname{depth}_{A}(a): = \min\{n :a\leq_{fin}r_n(A)\}
\end{equation}

\noindent \textbf{(A.3)}\ [{\bf Amalgamation}] Given $a$ and $A$
with $\operatorname{depth}_{A}(a)=n$, the following holds:
\begin{itemize}
		\item[(A.3.1)] $(\forall B\in [n,A])\ \ ([a,B]\neq\emptyset)$.
		\item[(A.3.2)] $(\forall B\in [a,A])\ \ (\exists A'\in [n,A])\ \
		([a,A']\subseteq [a,B])$.
\end{itemize}

\bigskip

\noindent \textbf{(A.4)}\ [{\bf Pigeonhole Principle}] Given $a$
and $A$ with $\operatorname{depth}_{A}(a) = n$, for every $\mathcal{O}\subseteq\mathcal{AR}_{|a|+1}$ there is $B\in
[n,A]$ such that $r_{|a|+1}[a,B]\subseteq\mathcal{O}$ or $r_{|a|+1}[a,B]\subseteq\mathcal{O}^c$.

\bteo[Todorcevic, \cite{todo}]\label{AbsEll}
	{\rm [{\bf Abstract Ellentuck Theorem}]}
	Any\linebreak $(\mathcal{R}, \leq, r)$ with $\mathcal{R}$ metrically closed and satisfying (A.1)-(A.4) is a
	topological\linebreak Ramsey space.
\eteo

\section{The topological Ramsey space $\mathcal P$}\label{def ramsey space polyhedra}
In this Section we will construct a new Ramsey space $\mc{P}$ and the set of its finite approximations
$\mc{AP}$.

\subsection{Definition of $\mc{P}$}\label{ssection def of P}
Consider pairs $(x,S_x)$ satisfying the following conditions:
\begin{itemize}
	\item[\Romnumero] $x\subseteq\mathbb N$,
	\item[\romnumero] $S_x\subseteq x^{[<\infty]}$ is {\bf hereditary}, i.e., $u\subseteq v\ \&\ v\in S_x \Rightarrow u\in S_x$, and
	\item[\romnumero] $\bigcup S_x = \bigcup\{u : u\in S_x\}= x$.
\end{itemize}
Given two such pairs $(x,S_x)$, $(y,S_y)$  let us define
\begin{equation}\label{eq preorder polyhedrons}
		(y,S_y) \leq (x,S_x)\ \Leftrightarrow \ y\subseteq x\ \&\ S_y\subseteq S_x.
\end{equation}
Let us write $\mc{AP}$ (resp. $\mc{P}$) for the set of all 
pairs $(x,S_x)$ satisfying properties (i), (ii), (iii) and such that $x$ is a finite (resp. an infinite) subset of $\N$. From now on, the elements of $\mc{P}$ will be written $(A,S_A), (B,S_B) \dots$, using capital letters. Let us define the preorder $\leq_{fin}$ and the partial order $\sqsubseteq$ on $\mathcal{AP}$ as  follows:
\begin{eqnarray}\label{eq preorder APN}
	(a,S_a) \leq_{fin} (b, S_b) & \Leftrightarrow & (a,S_a) \leq (b, S_b)\ \& \ \max(a) = \max (b)\\
	(a,S_a) \sqsubseteq (b, S_b) & \Leftrightarrow& a\sqsubseteq b\ \&\  (a,S_a) \leq (b, S_b)
\end{eqnarray}
Here we are using the same symbol $\sqsubseteq$ to indicate that the set 
$a$ is an initial segment of the set $b$. \\[2mm]
Given a pair $(A,S_A)\in\mc{P}$ and any subset $x\subseteq A$ (finite or countable),  let 
$S_A\upharpoonright x = \{u\cap x : u\in S_A\}$. 
In particular, if $n\in\N$ let 
$A\upharpoonright_n$ be the set of the first $n$ elements of $A$ and $S_{A\upharpoonright_n}=
S_A\upharpoonright (A\upharpoonright_n)$. 
The pair
\begin{equation}\label{eq initial segments polyhedra}
		r_n(A,S_A) = (A\upharpoonright_n, S_{A\upharpoonright_n})
\end{equation}
is the {\bf $n$th approximation} of $(A,S_A)$. Notice that
\begin{equation}\label{eq initial segments polyhedra-2}
		i\leq j\ \Rightarrow\ r_i(A,S_A) \leq r_j(A,S_A)\leq (A,S_A)\hskip1cm
		\forall i,j\in\N
\end{equation}
There is a well defined surjective function 
\begin{equation}\label{eq initial segments polyhedra 2}
		\mathcal P\times\mathbb N\ARROW{r} \mathcal {AP}\hskip1cm
		r((A,S_A), n) = r_n(A,S_A) 
\end{equation}

\subsection{$\mc{P}$ is a topological Ramsey space}\label{ssection P is Ramsey}
In the rest of this section we shall prove the following:
\bteo\label{polyhedra Ramsey space}
	$(\mathcal P, \leq, r)$ is a topological Ramsey space.
\eteo

The proof of Theorem \ref{polyhedra Ramsey space} will be divided into several lemmas, showing that $(\mathcal P, \leq, r)$ satisfies the conditions of 
the Abstract Ellentuck Theorem.
\blema\label{primeros tres}
	$(\mathcal P, \leq, r)$ satisfies axiom {\bf A.1}
	\begin{enumerate}
			\item For every $(A,S_A)\in\mathcal P$, $r_0(A,S_A) = \emptyset$.
			\item If $(A,S_A) \neq (B,S_B)$ then there exists $n$ such that $r_n(A,S_A) \neq r_n(B,S_B)$.
			\item If $r_n(A,S_A) = r_m(B,S_B)$ then $n=m$ and for every $i<n$, $r_i(A,S_A) = r_i(B,S_B)$.
	\end{enumerate}
\elema
\bdem
		Straightforward.
\edem

Hence each element of $\mathcal{P}$ can be identified with the sequence of its approximations. 
Next we consider $\mathcal P$ as a subset of the product space $\mathcal{AP}^{\mathbb N}$, 
regarding $\mathcal{AP}$ as a discrete space.

\blema\label{polyhedra closed space}
$\mathcal P$ is a closed subset of $\mathcal{AP}^{\mathbb N}$.
\elema
\bdem
		Consider the injection $\mc{P}\ARROW{\varphi}\mc{AP}^{\N}$ given by 
		\[
				\varphi(A,S_A)=(r_0(A,S_A),r_1(A,S_A),\dots).
		\]
		Let us show that $\varphi(\mc{P})$ is closed. Given a closure point 
		$\alpha= \left\{(a^j,S_{a^j})\right\}_{j\in\N}$ in $\overline{\varphi(\mc{P})}\subset\mc{AP}^{\N}$ 
		and a sequence $\left\{(A^k,S_{A^k})\right\}_{k\in\N}$ in $\mc{P}$,
		if $\left\{\varphi(A^k,S_{A^k})\right\}_{k\in\N}$ converges to $\alpha$ then 
		\[
					(\forall n\in\N)\ (\exists k\sub{n}\in\N)\hskip5mm
					k\geq k\sub{n}\Rightarrow (\forall j\leq n)\ r\sub{j}(A^k,S_{A^k})=(a^j,S_{a^j})
		\]
		Taking a strictly increasing sequence $k_n<k_{n+1}\forall n\in\N$ we get
		\[
					n>m\Rightarrow r\sub{m}(A^{k_n},S_{A^{k_n}})=(a^{k_m},S_{a^{k_m}})
		\] 
		Define $A=\underset{n\in\N}{\cup}\ a^{k_n}$ and $S_A=\underset{n\in\N}{\cup}\ S_{a^{k_n}}$.
		Then $(A,S_A)\in\mc{P}$ and $\varphi(A,S_A)=\alpha$ by construction.
\edem

The following two lemmas are straightforward; we leave the details to the reader.
\blema\label{cuarto axioma}
	$(\mathcal P, \leq, r)$ satisfies axiom {\bf A.2}
	\begin{enumerate}
		\item If $(A,S_A)\leq (B,S_B)$ then $\forall n\ \exists m$, $r_n(A,S_A)\leq_{fin} r_m(B,S_B)$.
		\item For every $(a,S_a)\in \mathcal{AP}$ the set $\{(b,S_b) : (b,S_b)\leq_{fin} (a,S_a)\}$ is finite.
		\item If $(a,S_a)\leq_{fin} (b,S_b)$ and $(c,S_c) \sqsubseteq (a, S_a)$ then there is $(d,S_d) \sqsubseteq (b, S_b)$ such that $(c,S_c) \leq_{fin}  (d, S_d)$.
	\end{enumerate}
\elema

Before stating Lemma \ref{quinto axioma} below, let us adapt from Section \ref{ssection Ramsey spaces} 
the definiton of basic open sets, for the \textbf{Ellentuck-like topology of $\mathcal P$}. 
These will be sets $[(a,S_a), (A,S_A)]$ such that $(B,S_B)\in [(a,S_a), (A,S_A)]$ iff 
\begin{equation}\label{eq Ellentuck topology of P}
      (B,S_B)\leq (A,S_A)\ \&\ (\exists n)\ r_n(B,S_B) = (a,S_a)
\end{equation}
In particular,
\begin{equation}\label{eq Ellentuck topology of P-2}
    [n, (A,S_A)] = [r_n (A,S_A), (A,S_A)]
\end{equation}
For $(a,S_a)\in\mathcal{AP}$ and $(A,S_A)\in\mathcal P$ such that $ [(a,S_a), (A,S_A)]\neq\emptyset$,
let us adapt the definition of {\bf depth} of $(a,S_a)$ in $(A,S_A)$ as follows  (cf. Axiom {\bf A.2}, p.4);
\begin{equation}\label{eq depth of a segment 2}
		\operatorname{depth}_{ (A,S_A)}(a,S_a): = \min\{n :(a,S_a)\leq_{fin}r_n(A,S_A)\}.
\end{equation}

\blema\label{quinto axioma}
	$(\mathcal P, \leq, r)$ satisfies axiom {\bf A.3}\\ 
	Let $n = \operatorname{depth}_{ (B,S_B)}(a,S_a)$.
	\begin{enumerate}
		\item  If $(A,S_A)\in [n, (B,S_B)]$ then $[(a,S_a), (A,S_A)]\neq\emptyset$.
		\item  For every $(A,S_A)\in [(a,S_a), (B,S_B)]$ there exists $(A',S_{A'})\in [n,(B,S_B)]$ such that $\emptyset\neq  [(a,S_a), (A',S_{A'})]\subseteq [(a,S_a), 		(A,S_A)]$ .
	\end{enumerate}
\elema
Given $n\in\N$ let 
\begin{equation}\label{eq def de APn}
      \mathcal{AP}_n := \{(a,S_a)\in\mathcal{AP} : |a| = n\}
\end{equation}
If $(a,S_a)\in\mathcal{AP}_n$ we say that the \textbf{length} of $(a,S_a)$ is $n$ or simply write 
$|(a,S_a)|=n$. Also, as in the general setting, for every  natural number $n$ write
\begin{equation}
r_n[(a,S_a), (A,S_A)] = \{r_n(B,S_B) : (B,S_B)\in [(a,S_a), (A,S_A)]\}
\end{equation}
Finally, we prove the following 
\blema\label{sexto axioma}
	Pigeonhole principle {\bf A.4} for $(\mathcal P, \leq, r)$:\\
	Let $n = \operatorname{depth}_{ (B,S_B)}(a,S_a)$, $k = |(a,S_a)|$ and 
	$c : \mathcal{AP}_{k+1} \rightarrow \{0,1\}$ be any partition. There exists
	$(A,S_A)\in [n,(B,S_B)]$ such that $c$ is constant in $r_{k+1} [(a,S_a), (A,S_A)]$.
\elema
\bdem
	Let 
	\[
			X = \{m\in B : m> max(a)\}.
	\] 
	For $i \in\{0,1\}$, let $$X_i = \{m\in X : c((a\cup\{m\}, 	S_B\upharpoonright a\cup\{m\})) = i\} .$$

	By the classical pigeonhole principle, there is $i_0 \in\{0,1\}$ such that $|X_{i_0}| = \infty$. So let 
	\[
			A = (B\upharpoonright n)\cup X_{i_0}\mbox{ and } S_A = S_B\upharpoonright A
	\] 
	Then $(A,S_A)\in [n,(B,S_B)]$ is as required.
\edem
Now we can prove that $(\mathcal P, \leq, r)$ is a topological Ramsey space:
\begin{proof}[Proof of Theorem \ref{polyhedra Ramsey space}]
In virtue of the abstract Ellentuck theorem, the required result follows from Lemmas \ref{primeros tres}, \ref{polyhedra closed space}, \ref{cuarto axioma}, \ref{quinto axioma}, \ref{sexto axioma}.
\end{proof}
\begin{rem}\label{obs Ellentuck's space as a subspace of P}\label{rem Ellentucks space as a subspace of P}
(\textbf{Ellentuck's space as a subspace of $\mathcal P$}) Notice that we can identify each $A\in\mathbb N^{[\infty]}$ with the pair $(A, A^{[\leq 1]})$. In this way, we can view $\mathbb N^{[\infty]}$ as a closed subspace of $\mathcal P$.
\end{rem}

 Recall the approximation function $i :\N\times\N^{[\infty]} \rightarrow \N^{[<\infty]}$, given by

$$i(n,A) = \ \mbox{the first } n\ \mbox{elements of } A.$$

Let $\mathcal E = (\N^{[\infty]}, \subseteq, i)$, where $\subseteq$ is the inclusion relation and $i$ is the approximation function defined above. For the space $\mathcal E $, the set of approximations is $\mathcal{AE} = \mathbb N^{[<\infty]}$. For every $a,b\in\mathbb N^{[<\infty]}$, \ $a\leq_{fin} b$ if and only if $a\subseteq b\ \&\ max(a)=max(b)$. Now we give an alternative proof to the well known fact that  $\mathcal E$ is a topological Ramsey space. 
\bcor\label{Ellentuck subspace} 
(Ellentuck \cite{ellen}, 1974) $\mathcal E = (\N^{[\infty]}, \subseteq, i)$ is a topological Ramsey space.
\ecor
\bdem
Fix $\mathcal X\subseteq \N^{[\infty]}$ with the Baire property with respect to the exponential topology of  $\mathcal E$. Since $\mathcal E$ is a closed subspace of $\mathcal P$, it is easy to show that the set 

\[
		\mathcal X' = \{(A,A^{[\leq1]}) : A\in \mathcal X\}\subset\mathcal P
\]

 has the Baire property with respect to the Ellentuck-like topology of $\mathcal P$
(\footnote{
  This can be deduced from two facts: \Letra Since $\mc{E}$ is closed in $\mc{P}$, every meager subset of 
  $\mc{E}$ is still meager in $\mc{P}$; and \letra Subsets of $\mc{P}$ with the Baire property form a 
  $\sigma$-algebra.
}). Given a nonempty neighborhood $[a,A]$ in $\mathcal E$, let $S_a = a^{[\leq 1]}$ and  $S_A = A^{[\leq 1]}$. Then, consider the neighborhood $[(a, S_a),(A, S_A)]$ in $\mathcal P$. Applying  Theorem \ref{polyhedra Ramsey space} we obtain $(B, S_B)\in[(a, S_a),(A, S_A)]$ such that 

$$[(a, S_a),(B, S_B)]\subseteq\mathcal X' \mbox{ or } [(a, S_a),(B, S_B)]\cap\mathcal X' = \emptyset.$$

\medskip

Notice that, by necessity, $S_B = B^{[\leq 1]}$. Hence,  $[a,B]\subseteq\mathcal X$ or $[a,B]\cap\mathcal X = \emptyset$.
\medskip

If $\mathcal X$ is meager with respect to the exponential topology of  $\mathcal E$ then the same argument works but in addition  the case $[a,B]\subseteq\mathcal X$ will never happen, by the meagerness of $\mathcal X$. This completes the proof.
\edem

From now on, we will refer to  $\mathcal E = (\N^{[\infty]}, \subseteq, i)$ as \textbf{Ellentuck's space}.

\subsection{Embeddings of ordered polyhedra}\label{ssection morphisms of simplexes}

 A \textit{finite ordered polyhedron} is a finite geometric polyhedron for which we have prefixed a linear order on the set of its vertices; it corresponds to a pair $(x,S_x)\in\mc{AP}$ considering $x$ with the natural order of $\mathbb N$. Hence,  $\mc{AP}$ can be understood as a subclass (in the sense of \cite{kpt}, for instace) of the class of finite ordered polyhedra. Similarly, $\mc{P}$ can be understood as a set of ordered polyhedra with a countable set of vertices. We call it  {\bf the Ramsey space of infinite countable ordered polyhedra}. The following will be useful to study this objects in relation to structural Ramsey theory.  An {\bf embedding} 
\[
    (x,S_x)\ARROW{f}(y,S_y)
\]
 is an injective function
$x\ARROW{f}y$ such that $u\in S_x\Rightarrow f(u)\in S_y$. It is a {\bf strong embedding} if  
$u\in S_x\Leftrightarrow f(u)\in S_y$.
A {\bf rigid embedding} is a strong embedding $(x,S_x)\ARROW{f}(y,S_y)$
such that $f$ is order-preserving: $i<j\Rightarrow f(i)<f(j)$. 
\blema\label{lema polyhedra embedded on simplexes}
	Each finite ordered polyhedron can be embedded in some $n$-simplex $\Delta$, for some $n\in\mathbb N$, and rigidly embedded in some subpolyhedron of $\Delta$.
\elema
\bdem
	If $(x,S_x)\in\mc{AP}$ is a finite ordered polyhedron let $n=|x|$, the cardinality of $x$, and write $x=\{x_0, x_1, \dots, x_{n-1}\}$ in increasing order. Define $x\ARROW{f}n$ by $f(x_j)=j$, for $j<n$. Then, $f$ induces an embedding
	$(x,S_x)\ARROW{f}\left(n,2\sup{n}\right)$
	which, geometrically,  is just the embedding of the polyhedron $K$ determined by 
	$(x,S_x)$ in the standard $n$-simplex $\Delta\sup{n}\subset\R{n+1}$. On the other hand, 
	$f$ also induces a rigid embedding $(x,S_x)\ARROW{f}\left(n,S_n\right)$
	where $S_n=\{u\subseteq n : \{x_j : j\in u\}\in S_x\}$.
\edem

\section{Finite polyhedra as a Ramsey class}\label{finite pol Ramsey}	
In this section we describe some basic concepts and results on Ramsey classes, Fra\"iss\'e theory and extremely 
amenability of automorphism groups, and prove that the class of finite ordered polyhedra is Ramsey. We prove that the automorphism group of its \Fraisse limit is extremely amenable, and state a universal property for ultrahomegeneous polyhedra. \\[2mm]
For the rest of this article, we will consider $L$-structures
$\mathbb{A} = \left\langle A, C\sup{\mathbb A}, R\sup{\mathbb A},F\sup{\mathbb A}\right\rangle$
on a fixed (first order) signature 
$L = \left\langle C, R,F\right\rangle$ of constants, relations, and fuctions symbols. 
Definitions such as morphisms, embeddings, isomorphisms, automorphisms, substructures, etc.,
can be found in the classical literature  
(for instance, see \cite{hodges}).

\subsection{Basic concepts}\label{ssection structures}\label{ssection substructures}
\label{ssection age of a structure}\label{ssection ultrahom structures}\label{ssection Fraisse classes}
The {\bf age} of an $L$-structure $\mathbb A$ is the class $\age{A}$ of all finite 
$L$-structures which are isomorphic to some substructure of $\mb{A}$. A structure $\mb{F}$ is \textbf{ultrahomogeneous} iff each isomorphism between any two finite substructures of $\mb{F}$
can be extended to some automorphism of $\mathbb F$. A {\bf \Fraisse structure} is a countable, locally finite,  
ultrahomogeneous structure.

\bteo[\Fraisse]\label{ultra son iso}
	Any two (infinite) countable  ultrahomogeneous $L$-structures having the same age are isomorphic.
\eteo

\bteo\label{Fraisse}
	 A non empty class of finite $L$-structures $\mc{C}$ is the age of a \Fraisse structure iff it satisfies:
	\begin{enumerate}
		\item $\mc{C}$ is closed under isomorphisms: If $\mb{A}\in\mc{C}$ and $\mb{A\cong B}$ then $\mb{B}\in\mc{C}$.
		\item  $\mc{C}$ is hereditary: If $\mb{A}\in\mc{C}$ and $\mb{B}\leq\mb{A}$ then $\mb{B}\in\mc{C}$.
		\item $\mc{C}$ contains structures with arbitrarily high finite cardinality.
		\item Joint embedding property: If $\mb{A,B}\in\mc{C}$ then there is $\mb{D}\in\mc{C}$ 
		such that $\mb{A\leq D}$ and $\mb{B\leq D}$.
		\item Amalgamation property: Given $\mb{A},\mb{B}_1,\mb{B}_2\in\mc{C}$ and embeddings 
		$\mb{A}\ARROW{f_i}\mb{B}_i$, $i\in\{1,2\}$, there is $\mb{D}\in\mathcal C$ and embeddings $\mb{B}_i\ARROW{g_i}\mb{D}$
		such that 
		$g_1\circ f_1 = g_2\circ f_2$.
	\end{enumerate}
	In such case, $\mathcal C$ is said to be a {\bf \Fraisse class}, and there exists a unique (up to isomorphism) countable \Fraisse structure $\mb{F}$ such that
	$\age{F}=\mc{C}$; this $\mb{F}$ is the {\bf \Fraisse limit} of $\mc{C}$ and we write $\mb{F}=\flim(\mc{C})$.
	
\eteo

\subsection{Ramsey classes of structures}\label{def clase ramsey}

Given $L$-structures $\mb{A,B,C}$ we write $\comb{B}{A}$ for the set of substructures of $\mb{B}$ which are isomorphic to $\mb{A}$. Given an integer $r>0$, if $\mb{A}\leq \mb{B}\leq\mb{C}$ then we write $\pigeonhole{A}{B}{C}{r}$ whenever for each $r$-coloring
	\[
			c :\comb{C}{A} \ARROW{} r
	\]
	of the set $\comb{C}{A}$, there exists $\mb{B'}\in \comb{C}{B}$ such that $\comb{B'}{A}$ is 
	monochromatic.
 
A \Fraisse class $\mathcal C$ has the {\bf Ramsey property} iff, for every integer $r>1$ and every $\mb{A,B}\in\mc{C}$ such that $\mb{A\leq B}$, there is $\mb{C}\in\mc{C}$ such that 
\[
			\ \  \ \pigeonhole{A}{B}{C}{r}\hskip1cm
\]

Also, remember that a topological group $G$ is \textbf{extremely amenable} or has 
\textbf{the fixed point on compacta property}, if for every continuous action of $G$ on a compact space
$X$ there exists $x\in X$ such that for every $g\in G$, $g\cdot x = x$. If $G$ is an extremely amenable 
group, then its \textit{universal minimal flow} is a singleton, a fact that is a remarkable result in 
Topological Dynamics. The following is an important characterization of the type of groups.
\bteo\label{aut fraise extrem edad ramsey}
	\cite{kpt}  Let $\mb{F}$ be a \Fraisse structure and $\mc{C}=\age{F}$. The polish group $\aut{F}$ is extremely amenable if and only if 
	$\mc{C}$ has the Ramsey property and all the structures of $\mc{C}$ are rigid.
\eteo

\subsection{Finite polyhedra as a Ramsey class}\label{def polyhedra ramsey category}

\medskip

Consider $L = <(R_i)_{i\in \mathbb N\setminus\{0\}}>$, a signature with an infinite number of relational 
symbols such that for each $i\in\mathbb N$ the arity of $R_i$ is $n(i) = i$. 

\medskip

A {\bf polyhedron} is a countable $L$-structure 
$\mb{A} = <A, (R_i^{\mb{A}})_{i\in \mathbb N\setminus\{0\}}>$ 
such that for each $\{a_1,\dots, a_i\}\subseteq A$, $(a_1,\dots, a_i)\in R_i$ if and only if 
$(a_{\sigma(1)},\dots, a_{\sigma(i)})\in R_i$, for every permutation $\sigma$ of the set $\{1, \dots, i\}$. 
Also, if $(a_1,\dots, a_i)\in R_i$, then for every $k\leq i$ and every subset $\{a_{j_1}, \dots, a_{j_k}\}$, 
we have $(a_{j_1}, \dots, a_{j_k})\in R_k$. Notice that if $\mathbb A$ is a finite $L$-structure, then there 
is a maximum arity $n = n(\mathbb A)$ such that 
$R_n^{\mathbb A}\neq\emptyset$ and $R_m^{\mathbb A} = \emptyset$, for every $m>n$. 
An \textit{ordered polyhedron} is a $L\cup\{<\}$-structure $\mathbb{A} = <A, (R_i^{\mathbb A})_{i\in \mathbb N\setminus\{0\}}, <^{\mathbb{A}}>$ such that $<A, (R_i^{\mathbb A})_{i\in \mathbb N\setminus\{0\}}>$ is a polyhedron and $<^{\mathbb{A}}$ is a total ordering on $A$.

\bigskip

Let $\mathcal{KP}_0$ be the class of finite polyhedra and $\mathcal{KP}$ the class of finite ordered\linebreak 
polyhedra.
It is easy to see that each pair $(x,S_x)\in\mathcal P\cup\mathcal{AP}$ 
is a countable $L\cup\{<\}$-structure whose universe is $x$ and in which $S_x$ is a countable 
family of relations over $x$. The notions of substructure, homorphism, etc, are induced by the embeddings defined in Section \ref{ssection morphisms of simplexes}. Furthermore, each one of these structures 
is rigid (ordered) by construction. 
\begin{rem}\label{hechos de los poliedros}  The following facts are straightforward:
    \begin{itemize}
	\item $\mathcal{AP}\subseteq \mathcal{KP}$. 
	\item For every $\mathbb A\in\mathcal{KP}$ there is $(a,S_a)\in\mathcal{AP}$ such that $\mathbb A\cong (a,S_a)$. Actually, $\mathcal{KP}$ is the closure of $\mathcal{AP}$ under isomorphisms.
    \end{itemize}
\end{rem}
We shall prove that the class $\mathcal{KP}$ is Ramsey in Theorem \ref{finite polyhedra are Ramsey} below. 
Before doing that we borrow the notation of \cite{nerod3}:
Let $\Delta = \{n_i\}_{i\in I}$ be a finite family of natural numbers. A \textit{set system of type} 
$\Delta$ is a structure $(X, \leq_X,\mathcal M)$ such that $(X, \leq_X)$ is a totally ordered set and 
$\mathcal M = \{\mathcal M_i\}_{i\in I}$ is such that $M\in X^{[n_i]}$, for every $M\in\mathcal M_i$.  
Given two set systems of type $\Delta$, $(X, \leq_X,\mathcal M)$ and $(Y, \leq_Y,\mathcal N)$ 
(with $\mathcal N = \{\mathcal N_i\}_{i\in I}$), we say that $(X, \leq_X,\mathcal M)$ is a \textit{subobject} 
of $(Y, \leq_Y,\mathcal N)$ whenever
\begin{itemize}
  \item $X\subseteq Y$, 
  \item $\leq_Y\upharpoonright X\times X\ =\ \leq_X$ and 
  \item $\mathcal M_i = \{M\in\mathcal N_i : M\subseteq X\}$.
\end{itemize}
Theorem A in \cite{nerod3} implies in particular that, for a fixed $\Delta$, the class of all sets systems 
of type $\Delta$ (together with all the embeddings) is Ramsey. It is easy to see that each 
$\mathbb A\in\mathcal{KP}$ is a set system of some type $\Delta_{\mathbb A}$. 
Actually, if $n=n(\mathbb A)$ is the maximum arity in $\mathbb A$ then 
$\Delta_{\mathbb A} = \{i \}_{i\leq n}$.
\begin{teo}\label{finite polyhedra are Ramsey}
    The class $\mathcal{KP}$ of all finite ordered polyhedra is Ramsey.
\end{teo}
\begin{proof}
  This follows as an application of Theorem A in \cite{nerod3}.\\[2mm] 
  Let $\mathbb A\leq \mathbb B\in\mathcal{KP}$ be given. Notice that
  $\Delta_{\mathbb A}$ is an initial segment of $\Delta_{\mathbb B}$, so we can assume that $\mathbb A$ and $\mathbb B$ have the same type (some relations in $\mathbb A$ can be empty). By Theorem A, there exists  a set system $\mathbb C = (X, \leq_X,\mathcal M)$,  $\mathcal M = \{\mathcal M_i\}_{i\in I}$, such that 
  \begin{itemize}
      \item $\Delta_{\mathbb C}=\Delta_{\mathbb B}$, and
      \item for every $r>1$, \ $\pigeonhole{A}{B}{C}{r}$. 
  \end{itemize}

Set $S_X = \{u\in\bigcup_{i\in I}\mathcal M_i : u \mbox{ is a face a copy of } \mathbb B \mbox{ inside } \mathbb C\}\cup X^{[\leq 1]}$ and let\linebreak $\mathbb D=(X, \leq_X, S_X)$. Then, $\mathbb D\in\mathcal{KP}$ and for every $r>1$, \ $\pigeonhole{A}{B}{D}{r}$. This completes the proof.
\end{proof}

\begin{rem}\label{finite polyhedra is age}
    The class of finite ordered polyhedra $\mc{KP}$ satisfies conditions (1),...,(5) of 
    Theorem \ref{Fraisse}; so it is the age of a \Fraisse structure. Let $\mb{P}= \flim(\mathcal{KP})$, the \Fraisse limit of $\mathcal{KP}$.
\end{rem}

\bcor

The automorphism group of $\flim(\mathcal{KP})$ is extremely amenable.

\ecor

\subsection{Geometric characterization of $\mb{P}= \flim(\mathcal{KP})$}
Now we will provide some arguments which are similar to those which 
arise in the construction of the \Fraisse limit of the class $\mathcal{KG}$ of finite graphs;
say $\Gamma= \flim(\mathcal{KG})$. There is a geometric characterization of $\Gamma$. 
For  each countable graph $\mb{G}=(V,E)$; we have $\mb{G}\cong\Gamma$  iff the following holds:
For any finite disjoint subsets of vertices $x,y\subset V$, there is some vertex 
$q\in V\setminus (x\cup y)$ such that $q$ is adjacent to all elements in $x$ 
and to none in $y$. See \cite[p.336-337]{hodges}. In order to show an analogous statement
for $\mb{P}$, we will start with two simple observations.\\[2mm]
Given a finite polyhedron $(a,S_a)$, we say that $T\subset\mc{P}(a)$ {\bf generates}
$S_a$  if $S_a=\{u:\exists v\in T(u\subseteq v)\}$. The $\subseteq$-minimal family generating $S_a$ is
$T_a=\max(S_a)$, the set of maximal subsets of $a$ in $S_a$ with respect to $\subseteq$.
Geometrically speaking, for any $T$ generating $S_a$, $T$ is a set of simplexes whose
amalgamation (union) in $\R{|a|+1}$ is the (geometric realization of the) polyhedron $(a,S_a)$; and $T_a$ is the family of maximal subsimplexes of $(a,S_a)$.\\[2mm]
A {\bf one-point extension} of a finite polyhedron $(a,S_a)$ is a finite
polyhedron $(b,S_b)$ such that $(a,S_a)\leq (b,S_b)$ and $b=a\cup\{p\}$ for some $p\notin a$. Then $p$ determines
a partition of $T_a$ into two classes: those $u\in T_a$ such that $u\cup\{p\}\in S_b$, and the other ones. 

\blema\label{lema induccion extensiones en un punto}
    For any countable polyhedron $\mb{A}=(A,S_A)$,  $\mb{A}$ is ultrahomogeneous iff the following condition holds:
\begin{itemize}
\item[{(*)}] For each finite polyhedron $(a,S_a)$, each embedding 
    $(a,S_a)\ARROW{f}\mb{A}$ and each one-point extension $(b,S_b)>(a,S_a)$;
    there exists an embedding\linebreak $(b,S_b)\ARROW{g}\mb{A}$ such that 
    $g\upharpoonright a=f$.
\end{itemize}
\elema

\bdem
    The direct implication is trivial. For the reciprocal, apply induction.
\edem

\bprop\label{prop pug poliedros}(Universal property for ultrahomogeneous polyhedra)

    A countable polyhedron $\mb{A}=(A,S_A)$ is ultrahomogeneous iff
    for any finite non empty disjoint subsets $x,y\subseteq S_A$, such that the elements
    of $x\cup y$ are not comparable by $\subseteq$, 
    there is some vertex $q\not\in \cup(x\cup y)$ of $\mb{A}$
    such that $u\cup\{q\}\in S_A$ $\forall u\in x$, and $u\cup\{q\}\not\in S_A$
    $\forall u\in y$.
\eprop
\bdem
    Fix $\mb{A}=(A,S_A)$. Let us show that the above geometric condition is\linebreak
    equivalent to Condition (*) in Lemma \ref{lema induccion extensiones en un punto}.\\
    $(\Rightarrow)$ Given a finite polyhedron $(a,S_a)$, 
    an embedding $(a,S_a)\ARROW{f}\mb{A}$, and a one-point extension 
    $(b,S_b)>(a,S_a)$ with $b=a\cup\{p\}$; consider the partition
    \[
	x_0=\{v\in T_a:v\cup\{p\}\in T_b\}\hskip1cm
	y_0=\{v\in T_a:v\cup\{p\}\not\in T_b\}  
    \]
    of $T_a$. Take $x=f(x_0)$, $y=f(y_0)$, so $x\cup y=f(T_a)$.
    By our assumption; there is $q\in A\backslash f(\cup T_a)$ such that
    $f(u)\cup\{q\}\in S_A$  $\forall u\in x_0$, and 
    $f(u)\cup\{q\}\not\in S_A$  $\forall u\in y_0$. Define $g(p)=q$.\\
    $(\Leftarrow)$ Let $x,y\subset S_A$ be as in the hypothesis. Take 
    $T=x\cup y$, $a=\cup T$ and $S_a$ the family generated by $T$. 
    Then $(a,S_a)$ is a finite polyhedron and $(a,S_a)<\mb{A}$. 
    Pick any $p\not\in a$. Let $b=a\cup\{p\}$, 
    $T_b=\{u\cup\{p\} :u\in x\}\cup y$, and $S_b$ the family generated by $T_b$.
    Then $(b,S_b)$ is a one-point extension of $(a,S_a)$; so there is an embedding
    $(b,S_b)\ARROW{g}\mb{A}$ satisfying $g\upharpoonright a\equiv id$. 
    Take $q=g(p)$. 
\edem

\section{The random polyhedron}\label{ssection infinite random polyhedron}
The article of P. Erdos \cite{erdos} is among the first approaches to the geometric properties of the 
random graph by means of probability methods. In this section we will study the universality property of countable random polyhedra. Here we will follow some arguments of \cite{rado}, where the property is also studied for random 
polyhedra.
\subsection{Definition of the random polyhedron}\label{ssection definition of the random polyhedron}
Hold a coin and assume that the probability of getting heads is $p=\frac{1}{2}$.  
Define a family $T\subseteq\N^{[<\infty]}$ as follows: For every 
$u\in \N^{[<\infty]}$ flip the coin, and say that $u\in T$ if and only if you get heads. Set 
\begin{equation}\label{eq caras del random poliedro}
      S := \{v : (\exists u\in T)\ v\subseteq u\}
\end{equation}
 and set $\omega : =\bigcup S$. 
Let us write, from now on, $S_{\omega}=S$ and $T_{\omega}=T$. So $(\omega, S_{\omega})$ is the amalgamation of the random family of simplexes $T_{\omega}$. 
\bteo\label{lema random polyhedron has as age the finite ones}
    Consider $(\omega, S_{\omega})$ as defined above. With probability 1, for each pair of finite, disjoint, non empty subsets $x,y\subset S_{\omega}$ 
    satisfying that the elements of $x\cup y$ are not comparable, there exists some     
    $q\in\omega\backslash\cup(x\cup y)$ such that $u\cup\{q\}\in S_{\omega}$ for all $u\in x$, and 
    $u\cup\{q\}\not\in S_{\omega}$  for all $u\in y$.    
\eteo
\bdem
  In  $2\sup{\N}$ take the product-metric topology, 
  the $\sigma$-algebra of the Borel sets $\mc{B}$ and the outer probability measure $P$ which extends, by 
  the Carath\'eodory method \cite{federer}, the probability of finite 
  coin flips. Given a finite subset $z\subset\N$ and   
  a finite tuple $a\in2\sup{z}$; the probability of the  
  basic open set $[a]=\left\{\chi\in2\sup{\N}:\chi\upharpoonright z=a\right\}$
  is $P([a])=\frac{1}{2\sup{|z|}}$. 
  Let $\N\sup{[<\omega]}\ARROW{X}2$ be any random process defining the family
  $T_{\omega}=X\sup{-1}(\{1\})$ which generates $(\omega,S_{\omega})$. Fix a bijection $\omega\sup{[<\infty]}\ARROW{\phi}\N$.
  Then $\chi=X\phi\sup{-1}$ is a random point in the probability space
  $\left(2\sup{\N},\mc{B},P\right)$. For any countable subset
  $A\subset2\sup{\N}$, we have $P(A)=0$. So, with probability 1, in a random countable flip one
  gets an infinite number of heads (see for instance \cite{feller}). 
  Since $T_\omega$ is infinite with probability 1, so is $\omega$.\\[2mm]
  Fix $x,y\subset S_\omega$ as in the hypothesis and let  $n=|x\cup y|$. For each 
  $q\in\omega\backslash\cup(x\cup y)$, consider 
  the proposition  
  \[
	  \varphi(q):=[(\forall u\in x)\ u\cup\{q\}\in S_\omega]\wedge[(\forall u\in y)\ u\cup\{q\}\not\in S_\omega].
  \]
  Let $z\sub{q}=\left\{\phi(u\cup\{q\}):u\in x\cup y\right\}$ and notice that
  the sets $z\sub{q}$, with $q$ ranging in $\omega\backslash\cup(x\cup y)$,
  are pairwise disjoint. Define
    \[
	a\sub{q}:z\sub{q}\ARROW{}2\hskip1cm
	a\sub{q}(\phi(u\cup\{q\}))=\left\{
	\begin{array}{lll}
	    1 && u\in x \\[2mm]
	    0 && u\in y
	\end{array}
	\right.
    \]
    Then the proposition $\varphi(q)$ is equivalent to the statement $\chi\upharpoonright z\sub{q}=a\sub{q}$. 
    The probability that $\varphi(q)$ holds is 
    $P(\chi\upharpoonright z\sub{q}=a\sub{q})=\frac{1}{2\sup{n}}$.
    Therefore, by the definition of $P$, for any finite subset 
    $F\subset[\omega\backslash\cup (x\cup y)]$, 
    \[
	P\left(\ \forall q\in F\ \neg\varphi(q)\ \right)=\left(1-\frac{1}{2\sup{n}}\right)\sup{|F|}
    \]
    So
    {\small\[
	\begin{array}{ll}
	 P(\ \forall q\in \omega\backslash \cup(x\cup y)\ \neg\varphi(q)\ ) 
	 & =\underset{q\not\in \cup(x\cup y)}{\prod} P(\ \chi\upharpoonright z\sub{q}\neq a\sub{q}\ )\\[4mm]
	 & \leq\underset{_{F\subset[\omega\backslash\cup (x\cup y)]}}{\text{LimInf}}
	 P\left(\ \forall q\in F\ \neg\varphi(q)\ \right)\\[4mm]
	 &=\underset{_{m\in\N}}{\text{LimInf}}\left(1-\frac{1}{2\sup{n}}\right)\sup{m}\\[4mm]
	 &=\Lim{m\rightarrow\infty}\left(1-\frac{1}{2\sup{n}}\right)\sup{m} =0. 
	\end{array}
    \]}

\medskip

    Therefore, $\exists q\in \omega\backslash \cup(x\cup y)\ \varphi(q)$ holds with probability 1. This completes the proof.
\edem
\bcor\label{unico random polyhedron}
   With probability 1, the following hold:
\begin{enumerate}
\item All infinite countable random polyhedra are ultrahomogeneous and isomorphic as countable structures.\label{unico random ultra}
\item $(\omega,S_\omega)$ contains
    finite simplexes of arbitrarily high dimension (cardinal).\label{random polyhedron contiene simplices grandes}
\item  Each finite polyhedron is rigidly embeddable on $(\omega,S_\omega)$. \label{finite polydra embed}
\end{enumerate}

\ecor
\bdem
   (1) By Theorem \ref{ultra son iso}, Proposition \ref{prop pug poliedros}  and Theorem \ref{lema random polyhedron has as age the finite ones}. (2) Take $y=\emptyset$, use Theorem \ref{lema random polyhedron has as age the finite ones} 
    and apply induction on $|x|$. (3)  By Part (\ref{random polyhedron contiene simplices grandes}) and Lemma \ref{lema polyhedra embedded on simplexes}.
\edem

Part (\ref{unico random ultra}) of Corollary \ref{unico random polyhedron} allows us to call $(\omega, S_{\omega})$ 
\textit{\underline{the} infinite countable random polyhedron}. We close this Section with the following.
\bteo\label{limite es el random poliedro}
    With probability 1,  $\mathbb P = \flim(\mathcal{KP})$ is an infinite ordered polyhedron which is isomorphic to  
    $(\omega, S_{\omega})$,  as a polyhedron, and to $(\mathbb Q,\leq)$, as an ordered set.
\eteo
\bdem
      By Corollary \ref{unico random polyhedron},
      $(\omega,S_\omega)$ is ultrahomogeneous and its age is $\mathcal{KP}$. 
      By Theorem \ref{ultra son iso}, $\mb{P}\cong(\omega,S_\omega)$ as a polyhedron.
      The second isomorphism is proved in a similar way, using the class of all finite
      linear orders.
\edem

\section{Topological Ramsey subspaces of $\mathcal P$}\label{subespacios}

In this section, for every integer $k>0$, we will define a topological Ramsey space $\mathcal P(k)$. 
It turns out that each $\mathcal P(k)$ will be a closed subspace of $\mathcal P$. 
In particular, $\mathcal P(1) = \mathcal E$, Ellentuck's space; and $\mathcal P(2)$ is a topological
Ramsey space whose elements are essentially the countably infinite ordered graphs. 
The corresponding set of approximations $\mathcal{AP}(2)$ is such that its closure under isomorphisms 
is essentially the class of finite ordered graphs, which is a Ramsey class (see \cite{nerod3}) whose \Fraisse 
limit is the ordered random graph. 
It is well-known that the automorphism group of the ordered random graph is, 
as in the case of the ordered random polyhedron, extremely amenable (see \cite{kpt}).

\subsection{The subspace $\mathcal P(k)$}
Given $k>0$, consider pairs of the form $(A,S_A)$ where:
\begin{itemize}
\item $A\in \mathbb N^{[\infty]}$.
\item $S_A\subseteq A^{[\leq k]}$.
\item $\bigcup S_A = A$.
\item $S_A$ is hereditary, i.e., $(u\subseteq v\ \& \ v\in S_A \Rightarrow u\in S_A$).
\end{itemize}
Let us define $\mathcal P(k)$ as the collection of all the pairs  $(A,S_A)$ as above. 
Consider the restrictions to $ \mathbb N\times\mathcal P(k)$ of the approximation function 
$r : \mathbb N\times\mathcal P \rightarrow \mathcal{AP}$ and the restrictions to  $\mathcal P(k)$  
and $\mathcal{AP}(k)$ of the pre-orders $\leq$ and $\leq_{fin}$ defined on $\mathcal P$  
and $\mathcal{AP}$.

\bteo\label{k-subspace ramsey}
For every integer $k>0$, the triplet $(\mathcal P(k), r, \leq)$ is a topological Ramsey space. In fact, it is a closed subspace of $(\mathcal P, r, \leq)$.
\eteo

\bdem
Given $k>0$, to show that  $(\mathcal P(k), r, \leq)$ is a  topological Ramsey space, proceed  as in the proof of Theorem \ref{polyhedra Ramsey space}. To show that is a closed subspace of $(\mathcal P, r, \leq)$, proceed  as in the proof of Corollary \ref{Ellentuck subspace}.
\edem

Actually, it is easy to show that given integers $k'>k>0$, $\mathcal P(k)$ is a closed subspace of $\mathcal P(k')$. As mentioned above,  $\mathcal P(1)$ is Ellentuck's space $\mathcal E$. 

\medskip

Now, let $\mathcal{KP}(k)$ denote the closure of $\mathcal{AP}(k)$ under isomorphisms. Proceeding as in the proofs of Theorems \ref{finite polyhedra are Ramsey} and \ref{finite polyhedra is age}, we obtain the following:

\begin{teo}\label{finite k-polyhedra are Ramsey}
For every $k>0$, the class $\mathcal{KP}(k)$ of all finite ordered $k$-polyhedra is Ramsey. Furthermore, $\mathcal{KP}(k)$ is the age of some \Fraisse structure.
\end{teo} 
\bcor

For every $k>0$, the automorphism group of $ \flim(\mathcal{KP}(k))$ is extremely amenable.

\ecor

\subsection{The random $k$-polyhedron}\label{random k-graph}

For every $k>0$  probabilistically define $(\omega, S^k_{\omega})$ with $S^k_{\omega}\subseteq\omega^{[\leq k]}$, proceeding just as in Section \ref{ssection infinite random polyhedron} (for defining $(\omega, S_{\omega})$). Notice that the corresponding versions of  Theorem \ref{limite es el random poliedro} and Corollary \ref{unico random polyhedron} now can be easily proved in this context. It turns out that the resulting pair $(\omega, S^k_{\omega})$ is characterized by the following, up to isomorphism:  with probability 1, $(\omega, S^k_{\omega})$ is isomorphic as a polyedron to $\flim(\mathcal{KP}(k))$. We then call $(\omega, S^k_{\omega})$ the \textbf{random $k$-polyhedron}. The following is the corresponding version of Theorem \ref{limite es el random poliedro}:

\bteo
Let $\mathbb P(k) = \flim(\mathcal{KP}(k))$, the Fra\"iss\'e limit of $\mathcal{KP}(k)$.
Then, with probability 1,  $\mathbb P(k)$ is isomorphic to $(\omega, S^k_{\omega})$,  as countable structure; and to $(\mathbb Q,\leq)$, as an ordered set.
\eteo

So $\mathbb P_k$ is the \textit{ordered} random $k$-polyhedron. It is also clear from the definitions that the random polyhedron contains an isomorphic copy of the random $k$-polyhedron.

\subsection{The ordered random graph}\label{random graph}

The case $k=2$ is of special notice. Observe that $\mathcal{KP}(2)$ is essentially the class of finite ordered graphs. Hence, as it is well-known, its Fra\"iss\'e limit $\mathbb P_2 = \flim(\mathcal{KP}(2))$ is the random ordered graph.

\section{Universal minimal flows}
Recall that every topological group $G$ has a universal minimal $G$-flow, unique up to isomorphism, which can be homomorphically mapped onto any other minimal $G$-flow. In this short Section we follow 
\cite{kpt} in order to calculate the universal minimal 
$G$-flow, when $G$ is the automorphism group of the random polyhedron or the automorphism group of the random $k$-polyhedron, $k\in\mathbb N\setminus\{0\}$.

\bdeff
Let $L$ be a signature with $\{<\}\subseteq L$, and put $L_0=L\setminus\{<\}$. Let $\mathcal K$ be a class of $L$-structures and put $\mathcal K_0=\mathcal K\upharpoonright L_0$. We say that $\mathcal K$ is \textbf{reasonable} if for every  $\mathbb A_0, \mathbb B_0\in\mathcal K_0$, every embedding $\pi : \mathbb A_0\rightarrow\mathbb B_0$, and every linear ordering $\prec$ on $A_0$ such that $\mathbb A = <\mathbb A_0,\prec\ >\ \in\mathcal K$, there is a linear ordering $\prec'$ on $B_0$, so that  $\mathbb B = <\mathbb B_0,\prec'\ >\in\ \mathcal K$ and  $\pi : \mathbb A\rightarrow\mathbb B$ is also and embedding. Also, we say that $\mathcal K$ satisfies the \textbf{ordering property} if for every  $\mathbb A_0\in\mathcal K_0$ there exists $\mathbb B_0\in\ \mathcal K_0$ such that for every linear ordering $\prec$ on $A_0$ and every linear ordering $\prec'$ on $B_0$, if $\mathbb A = <\mathbb A_0,\prec\ >\in\ \mathcal K$ and $\mathbb B = <\mathbb B_0,\prec'\ >\in\mathcal K$ then $\mathbb A\leq\mathbb B$.
\edeff

The class $\mathcal{KP}_0$ is the \textbf{reduct} of $\mathcal{KP}$; i.e. $\mathcal{KP}_0$ is the class of structures $\mathbb A_0$ 
obtained from structures $\mathbb A\in\mathcal{KP}$ by dropping the symbol $<^{\mathbb A}$. 
In a similar way,  for every $k\in\mathbb N\setminus\{0\}$ define the reduct $\mathcal{KP}_0(k)$ 
of $\mathcal{KP}(k)$. The classes $\mathcal{KP}_0$ and $\mathcal{KP}_0(k)$ are ages of Fra\"iss\'e structures. The non ordered random polyhedron $(\omega, S_{\omega})$ is the \Fraisse limit 
$\mathbb P_0$ of the class $\mathcal{KP}_0$, and the \Fraisse limit $\mathbb P_0(k)$ of the class $\mathcal{KP}_0(k)$ 
is the non ordered random $k$-polyhedron $(\omega, S^k_{\omega})$. In fact, if $L\supseteq\{\leq\}$ is the signature of the ordered polyhedra and $L_0=L\setminus\{\leq\}$ then, in the terminology of \cite{kpt}, $\flim(\mathcal{KP})\upharpoonright L_0=\flim(\mathcal{KP}_0)$  and $\flim(\mathcal{KP}(k))\upharpoonright L_0=\flim(\mathcal{KP}_0(k))$. These facts and Proposition 5.2 of 
\cite{kpt} imply the following:

\blema\label{KP reasonable}
    The classes $\mathcal{KP}$ and $\mathcal{KP}(k)$, $k\in\mathbb N\setminus\{0\}$, 
    are reasonable.\qed
\elema

Now, following \cite{kpt} again, in order to calculate the universal minimal flows for the groups $Aut(\mathbb P_0)$ and $Aut(\mathbb P_0(k))$ we need to show the following:

\blema\label{KP ordering}
    The classes $\mathcal{KP}$ and $\mathcal{KP}(k)$, $k\in\mathbb N\setminus\{0\}$, 
  satisfy the ordering property.
\elema
\bdem
We will prove that $\mathcal{KP}$  satisfies the ordering property. The rest can be done in an analogous way. We will proceed as in the proof of Theorem 2 of \cite{nerod2}, where it is proven that the class of all the finite ordered graphs satisfies the ordering property. So, fix $\mathbb A_0 = (a,S_a)\in\mathcal{KP}_0$ and a linear ordering $(a\leq)$. We will show that there exists $\mathbb B_0 = (b,S_b)\in\mathcal{KP}_0$ such that for every ordering $(b,\preceq)$ there exists a monotone mapping $f: (a,\leq)\rightarrow (b,\preceq)$ which is an embedding $(a,S_a)\rightarrow (b,S_b)$. Let $m$ be the number of finite polyhedra with set of vertices $a$ which are isomorphic to $(a,S_a)$. As is the case of finite graphs, if $m=1$ then we can take $(b,S_b)=(a,S_a)$, and therefore either $S_a=\emptyset$ or $S_a=a^{[\leq k]}$, for some $k\leq |a|$. So suppose $m>1$. By the Lemma in page 418 of \cite{nerod2}, for some large enough integer $n$, there exists a finite polyhedron $(c, S_c)$ such that its set of maximal 
faces $T_c$ satisfies $|T_c|=n^2$ and $(\forall u\in T_c)\, |u|=|a|$.   

Let $\mathcal F_c$ be the class of all finite polyhedra $(c, \hat{S}_c)$, with $c$ as set of vertices, such that:
\begin{enumerate}
\item For every $u\in T_c$, $(u, \hat{S}_c\upharpoonright u)\simeq (a, S_a)$.
\item $\hat{S}_c=\bigcup_{u\in T_c}\hat{S}_c\upharpoonright u$.
\end{enumerate}

Notice that $|\mathcal F_c|=m^{n^2}$. On the other hand, for every linear ordering $\preceq$ on $c$, there are less than $(m-1)^{n^2}+1$ elemets $(c, \hat{S}_c)$ of $\mathcal F_c$ for which there is no monotone embedding $$((a,\leq), S_a)\rightarrow ((c,\preceq), \hat{S}_c)).$$ Therefore, the cardinality of set of those $(c, \hat{S}_c)\in\mathcal F_c$ admitting an ordering $(c,\preceq)$ for which there exists no embeding $((a,\leq), S_a)\rightarrow ((c,\preceq), \hat{S}_c))$ is less than  $m!(m-1)^{n^2}$, which is $o(m^{n^2})$. This completes the proof.\edem

\medskip

Given a signature $L$ with $\{<\}\subseteq L$ and $L_0=L\setminus\{<\}$, let $\mathcal{K}$ be a reasonable class of $L$--structures  and  $\mathbb F = \flim(\mathcal{K})$. Let $F$ be the universe of  $\mathbb F$ and $\mathbb F_0=\mathbb F\upharpoonright L_0$. A linear ordering $\prec$ on $F$ is $\mathcal{K}$-\textbf{admissible} if  for every finite substructure $\mathbb A_0\leq\mathbb F_0$, we have $\mathbb A = <\mathbb A_0, \prec\upharpoonright A_0 >\ \in\mathcal K$, where $A_0$ is the universe of $\mathbb A_0$. By Lemma \ref{KP ordering}, in vitue of  Theorems \ref{finite polyhedra are Ramsey} and \ref{finite k-polyhedra are Ramsey} above, and Theorem 7.5(ii) of 
\cite{kpt} we obtain the following:

\bteo

Let $\mathbb P_0 = \flim(\mathcal{KP}_0) = (\omega, S_{\omega})$ be the random polyhedron and let $\mathbb P_0(k)= \flim(\mathcal{KP}_0(k)) = (\omega, S^k_{\omega})$, $k\in\mathbb N\setminus\{0\}$, be the random $k$-polyhedron. The following holds: 

\begin{enumerate}
\item The universal minimal $Aut(\mathbb P_0)$-flow is the metrizable $Aut(\mathbb P_0)$-flow  of \linebreak$\mathcal{KP}$-admissible orderings on $\omega$.
\item The universal minimal $Aut(\mathbb P_0(k))$-flow is the metrizable $Aut(\mathbb P_0(k))$-flow  of $\mathcal{KP}(k)$-admissible orderings on $\omega$.
\end{enumerate}
\eteo

\section{Final comment}	

The phenomena studied in this article reveal that, in general, there seems to be a tight relationship between a family of topological Ramsey spaces, and Ramsey classes of finite structures, extremely amenable automorphism groups,  universal minimal flows, etc. This raises several questions. For instance, consider the abstract setting introduced in  \cite{todo}. Given a Ramsey class $\mathcal K$ of (ordered) structures, what is the precise description of a topological Ramsey space  $\mathcal R$ (if any), such that $\mathcal K$ is the closure of $\mathcal{AR}$? Is it possible to characterize the family of topological Ramsey spaces $\mathcal R$ for which the class $\mathcal{AR}$ generates a Ramsey class of structures? On the other hand, a deeper study of the random polyhedron in itself and in relation to the random graph, from a wide point of view including approaches from model theory, graph theory, combinatorics, topology, dynamics and Ramsey theory is needed. We believe this work is a step in that direction.

\end{document}